\documentclass[final]{siamltex}
\usepackage{euscript,bbm,graphicx,epsfig,amsfonts,mathtools,epstopdf,amssymb}
\usepackage[ruled,lined,vlined,linesnumbered,algosection]{algorithm2e}
\usepackage{amssymb,latexsym}
\usepackage{psfrag}
 \usepackage[usenames]{color}
\usepackage[ruled,lined,vlined,linesnumbered,algosection]{algorithm2e}
\usepackage{url}
 
\def\cal{\mathcal}
\def\norm#1{\left\|#1\right\|}

\def\N{\mathbb{N}}
\def\R{\mathbb{R}}

\def\Rn{\R^n}

\def\Rnn{\R^{n\times n}}

\def\K{\mathbb{K}}

\def\BA{{\bf A}}

  \def\CH{{\cal H}}

  \def\CK{{\cal K}}

  \def\CP{{\cal P}}
  \def\CQ{{\cal Q}}
  \def\CR{{\cal R}}

\def\BAs{\BA{\kern-1.5pt}}

\def\CPs{\CP{\kern-0.8pt}}

\mathcode`\@="8000
{\catcode`\@=\active \gdef@{\mkern1mu}}

\def\mydate{\number\day\ {\ifcase\month \or January\or February\or
              March\or April\or May\or June\or July\or August\or
              September\or October\or November\or December\fi}
\number\year}
\def\vek#1{\mathbf{#1}}

\newcommand{\argmin}[1]{\underset{#1}{\text{{\rm argmin}}}}

\newcommand{\prn}[1]{\left(#1\right)}

\newcommand{\brac}[1]{\left[#1\right]}
\newcommand{\curl}[1]{\left\{#1\right\}}

\def\span{\text{span}}

\begin{document}

\title{A block MINRES algorithm based on the band Lanczos method \thanks{%
This version dated \today.}}

\author{Kirk M. Soodhalter\footnotemark[2]}

\renewcommand{\thefootnote}{\fnsymbol{footnote}}�
\footnotetext[2]{Industrial Mathematics Institute, Johannes Kepler University, 
Altenbergerstra{\ss}e 69, A-4040 Linz, Austria.
({\tt kirk.soodhalter@indmath.uni-linz.ac.at})}

\maketitle

\begin{abstract}
We develop a block minimum residual (MINRES) algorithm for symmetric indefinite matrices.  
This version is built upon the band Lanczos method that generates  one basis vector of the 
block Krylov subspace per iteration rather than a whole block as in the block Lanczos process.  
However, we modify the method such that 
the most expensive operations are still performed in a block fashion.  The benefit of
using the band Lanczos method is that one can detect breakdowns from
scalar values arising in the computation, allowing for a handling of breakdown which
is straightforward to implement.

We derive 
a progressive formulation of the MINRES method based on the band Lanczos process
and give some implementation details.  Specifically, a simple reordering of the steps allows us
to perform
many of the operations at the block level in order to take advantage of communication
efficiencies offered by the block Lanczos process. 
This is an important concern in the context of next-generation 
super computing applications.  

We also present a technique allowing us to maintain
the block size by replacing dependent Lanczos vectors with pregenerated
random vectors whose orthogonality against all Lanczos vectors is 
maintained.
Numerical results illustrate the performance on some sample problems.  
We
 present experiments that show how the relationship between right-hand
sides can effect the performance of the method.
\end{abstract}
\thispagestyle{empty}

\section{Introduction}\label{section.introduction}
We wish to efficiently solve 
\begin{equation}\label{eqn.AXB}
\vek A\vek X=\vek B
\end{equation}
where $\vek A\in\R^{n\times n}$ is a symmetric, indefinite matrix, and $\vek B\in\R^{n\times p}$
with $p$ right-hand sides.
  If $p=1$, Krylov subspace methods such as
the minimum residual MINRES method of Paige and Saunders \cite{Paige1975} have been shown to be 
effective.  For the case $p>1$, block Krylov subspace methods have been proposed; see, e.g.,  
\cite{FM.1997,L.2006,OLeary1980,Vital1990}.  In general, a block Krylov subspace method functions in
the same manner as a Krylov subspace method, but at each iteration the operator is applied
to a block of vectors rather than just one.  These methods generate $p$ new orthonormal basis 
vectors per iteration.
Many scalar operations become operations involving
small, dense matrices.  With these methods one can simultaneously solve the linear 
system for $p$ right-hand sides or solve a system with one right-hand side but over a block
Krylov subspace.  Though block methods increase the per iteration costs (as measured by
floating-point operation counts), they can be more efficient from the standpoint metrics related
to movement of data within the computer.

Our goals in this work are to develop a method which
\begin{enumerate}
	\item\label{item.goal1} solves (\ref{eqn.AXB}) over a block Krylov
subspace for any $p\geq 1$,
	\item\label{item.goal2} is designed to take advantage of the communication efficiencies of
block operations (when possible) but with greater ease of implementation,
	\item\label{item.goal3} is able to detect breakdowns through quantities arising in 
the computation,
	\item\label{item.goal4} and maintains the block size when a breakdown occurs.
\end{enumerate}
  Therefore, we seek a Lanczos-type method to generate the
block Krylov basis one vector per iteration, 
which is amenable to reordering of 
the steps to perform as many computations in blocks 
(e.g., sparse block operations or dense BLAS-3 operations)
as possible.  Furthermore, upon detection of a breakdown we prefer a
strategy which maintains the block size by replacing the dependent Lanczos
vector.

To this end we introduce a version of the MINRES algorithm for block Krylov subspaces
which satisfies our requirements.
This algorithm is built upon the band Lanczos process of Ruhe \cite{Ruhe1979}, 
which generates 
a basis for the block Krylov subspace one vector at a time rather than in a block fashion.  
Our algorithm can be considered a simplification of the algorithm presented in \cite{ABFH.2000}, 
which extends Ruhe's band Lanczos 
to generalize the nonsymmetric Lanczos process to the block setting in the case
that $\vek A$ is symmetric. However, we make modifications to execute some operations in a block fashion.  

We present the theoretical derivation needed to develop a minimum residual method based on
the band Lanczos procedure. 
We also discuss the simple modifications needed to execute some operations in a block fashion
 as well as some practical implementation details, to simplify the writing of the code.  
 To our knowledge this is the first paper to 
 provide implementation details of a block
minimum residual algorithm for symmetric matrices.\footnote{Matlab implementation available at \url{http://math.soodhalter.com/software.php}}

In the next section, we introduce notation and give a brief review of Krylov subspace methods (both non-block
and block).  In Section \ref{section.bminres},
we derive a version of the block minimum residual method built upon
the band Lanczos method.
In Section \ref{section.block-minres-impl}, we derive the progressive formulation of this 
method in detail, which is built to take advantage of the potential
memory savings afforded by the method. In Section \ref{section.dep-rhs}, 
we present modifications to our implementation which accommodate the occurrence of
exact or inexact dependence of a candidate block Krylov subspace basis vector.
In Subsection \ref{subsection.dyn-sub}, we develop a technique to 
maintain the block size when breakdown occurs.
In Section \ref{section.conv}, we discuss convergence properties
of block methods.
In Section \ref{section.impl}, special attention is given to how certain data 
is stored to keep the scheme as simple as possible.  
 In Section \ref{section.num}, we present numerical results.

\section{Preliminaries}\label{section.prelim}

For the matrix $\vek A$
and starting vector $\vek r_{0}$ recall that we generate an orthonormal basis for 
the Krylov subspace
\[
\CK_{j}(\vek A,\vek r_{0}) = \text{span}\curl{\vek r_{0}, \vek A\vek r_{0},\ldots, \vek A^{j-1}\vek r_{0}}
\]
with the Arnoldi process.  Let $\vek V_{j}\in\R^{n\times j}$ be the 
matrix with orthonormal columns consisting of this basis.  
Then we have the Arnoldi relation
\begin{equation}\label{eqn.arnoldi-relation}
	\vek A\vek V_{j} = \vek V_{j+1}\overline{\vek H}_{j}
\end{equation}
with $\overline{\vek H}_{j}\in\R^{(j+1)\times j}$; see, e.g., 
\cite[Section 6.3]{Saad.Iter.Meth.Sparse.2003}
and \cite{szyld.simoncini.survey.2007}.
In the case $p=1$, we can solve (\ref{eqn.AXB}) with a Krylov subspace iterative method.  Suppose
$\vek x_{0}$ is an initial approximation, and $\vek r_{0}=\vek B-\vek A\vek x_{0}$ is the initial 
residual.  At iteration $j$, we can compute the minimum residual correction 
$\vek t_{j}\in\CK_{j}(\vek A,\vek r_{0})$ satisfying
$\vek t_{j} =\argmin{\vek t\in\CK_{j}(\vek A,\vek r_{0})}\norm{\vek B-\vek A(\vek x_{0}+\vek t)}$
by solving the equivalent small least squares problem
$\vek y_{j}=\argmin{\vek y\in\R^{j+1}}\norm{\norm{\vek r_{0}}\vek e_{1}^{(j+1)} - \overline{\vek H}_{j}\vek y}$
and setting $\vek t_{j} = \vek V_{j}\vek y_{j}$.
Implementations of minimum residual methods include GMRES \cite{Saad.GMRES.1986} in the 
nonsymmetric case and MINRES \cite{Paige1975} in the symmetric case.

Much has been written about the solution of linear systems with multiple right-hand sides. 
Extending the framework of a Krylov method to the block right-hand side 
setting involves generalizing the machinery to deal with block vectors; see, 
e.g., \cite[Page 208]{Saad.Iter.Meth.Sparse.2003}.  

Let $\vek V_{1}\in\R^{n\times p}$ be a 
matrix with orthonormal columns.
At step $j$, the block Arnoldi process
generates an orthonormal basis for a block 
Krylov subspace 

\[
\K_{j}(\vek A,\vek V_{1}) = \CK_{j}(\vek A, \vek v_{1}^{(1)})+ \CK_{j}(\vek A, \vek v_{1}^{(2)}) + \cdots + \CK_{j}(\vek A, \vek v_{1}^{(p)}).
\]
where 
$\vek V_{1} = \brac{\begin{matrix} \vek v_{1}^{(1)}, &\vek v_{1}^{(2)}, &\ldots ,& \vek v_{1}^{(p)} \end{matrix}}$. 
In this setting sparse matrices act on a block of vectors per iteration.
For this discussion, we assume for now that 
$\dim \K_{m}(\vek A,\vek V_{1}) = mp$,
i.e., no 
linear dependent block Arnoldi vectors are generated.  We will return to the case of dependence later.

A strength of generalizing the Arnoldi process (or Lanczos) is that the block level
operations (e.g., BLAS-3 operations for dense matrices and block operations
for sparse matrices) have been shown to be quite efficient, when 
 measured in metrics relevant in a high performance computing environment, i.e., amount of 
data moved through memory, frequency of cache missing, and the number of floating
point operations performed on a unit of data while it is in cache.
These considerations have led to the broader goal to design communication avoiding Krylov subspace
methods; see, e.g.,  \cite{H.2010}.  In next-generation supercomputing machines the movement of data within the machine 
(e.g., from main memory onto the cache) increasingly will represent the dominant computational cost, 
and algorithms should be 
judged according to an appropriate data movement metric \cite{DS.2012}.  When judged according to 
such metrics the per iteration data movement costs of a block method are only marginally more expensive 
than their single vector counterparts (both dense, BLAS-3 operations and sparse block operations).  
Therefore, we minimize the residual over a larger constraint space without a concomitant increase 
in computational costs (related to the movement of data).  For further work on this topic, see, e.g., 
\cite{BDJ.2006,PSS.2013}.

One can also generate the block Krylov subspace one vector at a time
using the band Lanczos process proposed by Ruhe 
\cite{Ruhe1979}. 
At each iteration one 
matrix-single-vector product is performed as opposed to a 
matrix-block-vector product in a block-level method.  
It proceeds in a similar fashion to the single-vector Arnoldi process but starts with $p$ vectors
against which the new vector must be orthogonalized instead of one.
We derive our algorithm from this process
but with certain operations performed in a block fashion.  
We  must adopt a notation which is compatible with the
single vector per iteration nature of the band Lanczos process.
Thus the initial block of normalized vectors called $\vek V_{1}$ before
is renamed 
$\vek U_{p}=\brac{\begin{matrix} \vek u_{1}, \ldots , \vek u_{p} \end{matrix}} $, denoting that we start with the first $p$
orthonormal vectors.

Beginning with no symmetry assumption on $\vek A$, the 
band Arnoldi process (see, e.g., \cite{F.1995,FM.1997}) performs the same orthogonalization as the 
block method, only one vector at a time.  We denote the matrix with the 
first $j$ band Arnoldi vectors as columns $\vek U_{j}\in\R^{n\times j}$ where for $j < p$,
$\vek U_{j}$ has only the first $j$ starting vectors as columns.

This 
algorithm allows one to detect a breakdown from the scalar
quantities generated by the band Lanczos process. By reordering
the computations the band Lanczos algorithm can be formulated with many of the same 
block level operations as the block Lanczos algorithm, e.g., the operator is applied to a block of vectors
every $p$ iterations while maintaining the ease with which we detect breakdown in the 
band Lanczos algorithm.

To describe the Arnoldi relation in this setting, we must take care as the iteration number $j$ does not
match the dimension of the block Krylov subspace.  At iteration $j$, we generate the $(j+p)$th band
Arnoldi vector.  
At this iteration we have the band Arnoldi relation
\begin{equation}\label{eqn.Ruhe-bl-Arnoldi-rel}
\vek A\vek U_{j} = \vek U_{j+p}\overline{\vek H}_{j}.
\end{equation}
The banded Hessenberg matrix $\overline{\vek H}_{j}\in\R^{(j+p)\times j}$ has $p$ lower subdiagonal entries per column 
and has the structure
\[
	\overline{\vek H}_{j} = \brac{\begin{matrix}\vek H_{j} \\ \vek H_{p\times j}\end{matrix}},
\]
where $\vek H_{j}$ is a square $j\times j$ matrix 
satisfying the identity
\begin{equation}\label{eqn.sqr-hess-ident}
\vek H_{j} = \vek U_{j}^{\ast}\vek A\vek U_{j}.
\end{equation} 
Observe that $\vek H_{p\times j}$ only has nonzero entries in 
the last $p$ columns with structure
$\vek H_{p\times j} = \brac{\begin{matrix}\vek 0_{p\times (j-p)} & \vek C_{j}\end{matrix}}$
where $\vek C_{j}\in\R^{p\times p}$ is upper triangular. 
At iteration $j$ the dimension of the subspace built  is $j+p$.

To unambiguously describe the subspace at each iteration, we identify it with
the pair $(k,m)$ determined uniquely by
$j=(k-1)p+m$ with $0\leq m<p$.
As shown in \eqref{eqn.block-Krylov-subspace}, 
this pair is used to describe the 
block Krylov subspace built by the band Lanczos process.
The subspace that has been generated at iteration $j$ is the sum of 
Krylov subspaces generated by each column of $\vek U_{p}$, i.e.,
\begin{equation}\label{eqn.block-Krylov-subspace}
\K_{k,m}(\vek A,\vek U_{p}) = \underbrace{\sum_{\ell=1}^{m} \CK_{k+1}(\vek A,\vek u_{\ell})}_{\varnothing\mbox{\ when\ }m=0}+ \sum_{\ell=m+1}^{p}\CK_{k}(\vek A,\vek u_{\ell}) = \CR\prn{\vek U_{j}}.
\end{equation} 
Initially we have the identity, 
$\K_{1}(\vek A,\vek U_{p}) = \K_{1,0}(\vek A,\vek U_{p})$.
After $p$ iterations we have the following sequence of 
nested subspaces,
\begin{equation}\nonumber
	\K_{1,1}(\vek A,\vek U_{p})\subseteq \K_{1,2}(\vek A,\vek U_{p})\subseteq \cdots\subseteq \K_{1,p}(\vek A,\vek U_{p}).
\end{equation}
We observe that, in fact, 
$\K_{1,p}(\vek A,\vek U_{p}) = \K_{2,0}(\vek A,\vek U_{p})$,
since $m<p$ must hold.
The Krylov subspaces in the sum
(\ref{eqn.block-Krylov-subspace}) for the first $m$ right-hand sides 
$\curl{\vek u_{1},\cdots,\vek u_{m}}$
are of dimension $k+1$
and the remaining are of dimension $k$.
For $j=(K-1)p$, a multiple of the block size, 
the band Lanczos process has produced an orthonormal basis 
spanning the $K$th {\em block} Krylov subspace generated by $\vek A$ and $\vek U_{p}$, i.e., 
\begin{equation}\label{eqn.block-band-equiv}
	\K_{K}(\vek A,\vek U_{p}) = \K_{K ,0}(\vek A,\vek U_{p}).  
\end{equation}
At each iteration one of the subspaces in the 
sum (\ref{eqn.block-Krylov-subspace}) increases by one dimension.   

Similar to the symmetric Lanczos relation in the case of a single-vector Krylov method, observe that 
if $\vek A$ is symmetric the relation (\ref{eqn.sqr-hess-ident})
implies that ${\vek H}_{j}$ is also symmetric.  Due to the banded Hessenberg structure of ${\vek H}_{j}$, 
we see that  $\overline{\vek H}_{j}$ is a banded matrix with 
$p$ superdiagonal entries and $p$ subdiagonal entries per column.  
This structure implies that the orthogonalization process requires
only the 
most recent $2p$ basis vectors in 
order to compute $\vek u_{j+p}$.  We have the $2p+1$ term recurrence relation 
\begin{equation}\label{eqn.block-Lanczos-recurr}
	\vek A\vek u_{j} = \sum_{\ell=\min\curl{1,j-p}}^{j+p}h_{\ell,j}\vek u_{\ell}.
\end{equation}

Due to symmetry we do not need to compute $h_{\ell,j}$ where $\ell<j$ since it was computed previously 
as ${h_{j,\ell}}$.  
 This yields Algorithm \ref{alg.block-Arnoldi-Ruhe}, Ruhe's band 
Lanczos method.

\begin{algorithm}
\caption{The band Lanczos Process}
\label{alg.block-Arnoldi-Ruhe}

\SetKwInOut{Input}{Input}\SetKwInOut{Output}{Output}
\Input{$\vek A\in \Rnn$ symmetric, $\vek U_{p}\in\R^{n\times p}$, $\vek U_{p}^{\ast}\vek U_{p} = \vek I_{p}$}
\Output{$\vek U_{j+p}\in\R^{n\times (j+p)}$, $\vek U_{j+p}^{\ast}\vek U_{j+p} = \vek I_{j+p}$ and 
$\overline{\vek H}_{j}\in\R^{(j+p)\times m}$,  $\overline{\vek H}_{j}$ is symmetric 
with $p$ lower subdiagonal entries}

\For{$\ell=1, 2,\ldots, j$}{
	Compute $\vek w:=\vek A \vek u_{\ell}$\\
		\For{$i=\max\curl{1,p-\ell},\ldots,p+\ell-1$}{
			$h_{i,\ell}:=\vek u_{i}^{\ast}\vek w$\\
			$\vek w \leftarrow \vek w - h_{i,\ell}\vek u_{i}$\\
		}
	Compute $h_{p+\ell,\ell}:=\norm{\vek w}_{2}$ and $\vek u_{p+\ell}:=\vek w/h_{p+\ell,\ell}$
}

\end{algorithm}

It should be noted; our aim is in contrast to the goals stated in the dissertation of Loher \cite{L.2006} 
in which
the author extended the work of Aliaga et al \cite{ABFH.2000} to a fully block
nonsymmetric Lanczos-based method, preferring the flexibility offered by
a block method, e.g., with regard to look-ahead and deflation.  Furthermore, our approach 
can be considered as an alternative to the fully block approach of O'Leary \cite{OLeary1980}.
Schmelzer analyzed fully block MINRES and SYMMLQ in \cite{S.2004}.
The strategy advocated in the present work was commented upon in \cite{L.2006} as 
an alternative strategy one could pursue.
The flexibility of the fully block methods with regard to breakdowns come at the price of a more
complicated implementation.  Here, we sacrifice some of this flexibility in exchange for some simplicity
of implementation.  

We end by describing some nomenclature and notation. 
We call a vector with multiple columns, such as $\vek B$ when $p>1$, a \textit{block vector}.    
Boldface upper-case letters are used to denote matrices, including block vectors.  Boldface 
lower-case letters will denote column vectors.  We denote the Euclidean norm by $\norm{\cdot}$.  
For a square, nonsingular matrix $\vek A$, we will denote 
the condition number associated with the $2$-norm $\kappa(\vek A) = \norm{\vek A}\norm{\vek A^{-1}}$.  
When identifying an equation as a QR-factorization, we will use the convention that the right-hand side 
of the equation is the QR-factorization of the left-hand side of the equation.  We denote the $k\times k$ 
identity matrix $\vek I_{k}$.  We also use the Matlab indexing notation to indicate
a range of rows or columns of a matrix, e.g., $\vek M{(i:j,:)}$ is the submatrix
containing rows $i$ to $j$ and all columns of $\vek M$.
We have similarly for a product of matrices $\prn{\vek M\vek N}_{i:j,:}$
to avoid ambiguity.  For a matrix $\vek M$, we denote its range
 (i.e., the span of the columns) by $\CR(\vek M)$.

Since the word {\em deflation} has more than one meaning in our community we will refer to 
the process of removing dependent vectors to maintain a 
linearly independent basis in a block Krylov subspace method 
simply as \textit{removal of dependent vectors}.

\section{A Block Minimum Residual Method}\label{section.bminres}
We derive a minimum residual algorithm based on the band Lanczos process.  
If we begin with an initial guess $\vek X_{0}$,
at the $j$th step the following method will produce an approximation $\vek X_{j}\in\R^{n\times p} $ 
such that for each $0<i\leq p$, the residual $\norm{\vek b^{(i)} 
- \vek A\vek x_{j}^{(i)}}$ is minimized over the space $\K_{k,m}(\vek A,\vek F_{0})$, where $\vek x_{j}^{(i)}$ 
is the $i$th column of $\vek X_{j}$, and $\vek F_{0}=\vek B-\vek A\vek X_{0}$ is the 
initial residual. 

At step $j$, we minimize each column of the block residual
$\vek F_{j}=\vek B-\vek A\vek X_{j}$
over $\K_{k,m}(\vek A,\vek F_{0})$.  
Following the 
development of MINRES presented in \cite{Greenbaum1997a} we 
can derive a block MINRES algorithm based on the band Lanczos process. 
Let  $\vek E^{(j)}_{1}\in \R^{(j+p)\times p}$ be the matrix containing the first $p$ 
columns of $\vek I_{j+p}$.  Observe that 
\begin{equation}\label{eqn.E1.block}
\vek E_{1}^{(j)} = \brac{\begin{matrix}\vek E_{1}^{(j-1)}\\ \vek 0_{1\times p}\end{matrix}}.
\end{equation}  
Given $\vek F_{0}$ we can normalize it by computing the 
economized QR factorization 
\begin{equation}\label{eqn.init-QR}
\vek F_{0} = \vek U_{p}\vek S,
\end{equation}
where $\vek U_{p}\in\R^{n\times p}$ has
orthonormal columns and $\vek S\in\R^{p\times p}$ is upper triangular.

At step $j$ of band Lanczos process, we have the 
QR factorization $\overline{\vek H}_{j} = \vek Q_{j}\overline{\vek R}_{j}$ such that 
$\vek Q_{j}\in\R^{(j+p)\times (j+p)}$ is unitary, and $\overline{\vek R}_{j}\in
\R^{(j+p)\times j}$ is upper triangular.  The matrix $\overline{\vek R}_{j}$ has a simple block structure,
\[
\overline{\vek R}_{j} = \brac{\begin{matrix}\vek R_{j}\\ \vek 0_{p\times j}\end{matrix}},
\]
where $\vek R_{j}$ is a square, upper triangular, $j\times j$ matrix.
Let $\vek f_{j}^{(i)}$ be the $i$th column of $\vek F_{j}$, the $j$th block residual.
The minimization of $\norm{\vek f_{j}^{(i)}}$ can be rewritten as
\begin{eqnarray}
\label{eqn.min-resid} \norm{\vek f_{j}^{(i)}} & = & \min_{\vek x\in \vek x_{0}^{(i)} + \K_{k,m}(\vek A,\vek B)}\norm{\vek b^{(i)} - \vek A\vek x}\nonumber\\
                         & = & \min_{\vek y \in\R^{j}}\norm{\vek f_{0}^{(i)} - \vek A\vek U_{j}\vek y}\nonumber\\
                         & = & \min_{\vek y \in\R^{j}}\norm{\vek U_{p}\vek S\vek e^{(i)}_{p} - \vek A\vek U_{j}\vek y}\nonumber\\
                         & = & \min_{\vek y \in\R^{j}}\norm{\vek U_{j+p}\vek E^{(j)}_{1}\vek S\vek e^{(i)}_{p} - \vek U_{j+p}\overline{\vek H}_{j}\vek y}\nonumber\\
                         & = & \min_{\vek y \in\R^{j}}\norm{\vek E^{(j)}_{1}\vek S\vek e^{(i)}_{p} - \overline{\vek H}_{j}\vek y}\nonumber\\
                         & = & \min_{\vek y \in\R^{j}}\norm{\vek Q_{j}^{\ast}\vek E^{(j)}_{1}\vek S\vek e^{(i)}_{p} - \overline{\vek R}_{j}\vek y}.\label{eqn.resid-col-min}
\end{eqnarray}
We remind the reader that the upper triangular matrix $\vek S$ coming from 
(\ref{eqn.init-QR}) serves the same role as the norm of the initial
residual in single-vector Krylov methods. 

We can solve the normal equations individually for each right-hand side, or we can solve for all 
right-hand sides simultaneously, i.e., 
\begin{equation}\nonumber
	\vek Y_{j} = \vek R_{j}^{-1}(\vek Q_{j}^{\ast}\vek E^{(j)}_{1}\vek S)_{1:j,:}\mbox{\ \ with\ \ }\vek X_{j} = \vek X_{0} + \vek U_{j}\vek Y_{j}
\end{equation} 
Similar to the development of MINRES for one right-hand side in \cite{Paige1975} we define
$\overline{\vek Z}_{j} = \vek Q_{j}^{\ast}\vek E^{(j)}_{1}\vek S$.
The first $j$ rows of $\overline{\vek Z}_{j}$ define the coefficients of the correction in the basis of search directions 
defined by 
$\vek M_{j} = \vek U_{j}\vek R_{j}^{-1}$.
Observe that the columns of $\vek M_{j}$ successively span the same subspaces as the columns of 
$\vek U_{j}$ due to the upper triangular
structure of $\vek R_{j}^{-1}$. We denote the block vector of search direction 
coordinates 
$\vek Z_{j} = \overline{\vek Z}_{j}\prn{1:j,:}$.
The block minimum residual approximation at step $j$ is 
\begin{eqnarray}
	\vek X_{j} & = & \vek X_{0} + \vek U_{j}\vek Y_{j}\nonumber\\
	                    & = & \vek X_{0} + \vek U_{j}\vek R_{j}^{-1}(\vek Q_{j}^{\ast}\vek E^{(j)}_{1}\vek S)_{1:j}\nonumber\\
	                    & = & \vek X_{0} + \vek M_{j}\vek Z_{j}\label{eqn.minres-approx}.
\end{eqnarray}
It remains to show that, as in the case of MINRES, this indeed leads to a progressive
formulation.  As in the single right-hand side case, 
a computed residual (also sometimes called the \textit{recursive 
residual}) is available,
\begin{equation}\label{eqn.comp-resid}
	\norm{\vek f_{j}^{(i)}} = \norm{\overline{\vek z}_{j}^{(i)}\prn{j+1:j+p}},
\end{equation}
where $\overline{\vek z}_{j}^{(i)}$ is the $i$th column of $\overline{\vek Z}_{j}^{(i)}$.
This can be derived from (\ref{eqn.resid-col-min}), which can be rewritten as 
\begin{equation}\label{eqn.resid-col-min-cpt}
	\norm{\vek f_{j}^{(i)}}  = \min_{\vek y \in\R^{j}}\norm{\overline{\vek z}_{j}^{(i)}- \overline{\vek R}_{j}\vek y}.
\end{equation}
As we assume here that there has been no breakdown in the band Lanczos process, $\vek R_{j}$
is nonsingular.  Thus, (\ref{eqn.resid-col-min-cpt}) can be satisfied exactly in the first $j$ rows.  Due to the 
structure of $\overline{\vek R}_{j}$, we have that the residual is simply the norm of the last $p$ entries of 
$\overline{\vek z}_{j}^{(i)}$, i.e., (\ref{eqn.comp-resid}).

\section{Block MINRES for Symmetric Linear Systems}\label{section.block-minres-impl}
To obtain a storage-efficient block MINRES algorithm based on the band Lanczos method we must discuss the structure of $\vek R_{j}$.  This matrix is the 
upper $j\times j$ block of $\overline{\vek R}_{j}$ which is obtained from the QR-factorization of 
$\overline{\vek H}_{j}$.  As the lower subdiagonal of $\overline{\vek H}_{j}$ has $p$ nonzero entries, 
we obtain this factorization 
using Householder reflections
To each new column of $\overline{\vek H}_{j}$, we must apply all previous
reflections.  This procedure adds 
 to the new column at most $p$ new nonzero superdiagonal entries.
 As a result the upper triangular $\vek R_{j}$ has at most $2p$ superdiagonal entries per column.
   
The identity $\vek M_{j}\vek R_{j} = \vek U_{j}$ yields
the relationship between the band Lanczos 
vectors and the search directions,
\begin{eqnarray}
r_{m,j}\vek m_{m} + r_{m+1,j}\vek m_{m+1} + \cdots + r_{j,j}\vek m_{j} & = & \vek u_{j},\label{eqn.search-dir-comp} 
\end{eqnarray}
where $m= \min\curl{1,j-2p}$.
Thus to compute $\vek m_{j}$ we need $\vek u_{j}$ and the $2p$ previous search directions.  

The Householder reflections must also be applied to $\vek E_{1}^{(j)}\vek S$
to construct the residual according to \eqref{eqn.resid-col-min}.
Let $\CH_{i}^{(\ell)}\in\R^{\ell\times\ell}$ be the Householder reflection annihilating
the entries in the 
$i$th subdiagonal of $\overline{\vek H}_{\ell}$.
From (\ref{eqn.E1.block}) we have that $\vek E_{1}^{(j-1)}\vek S$ is a submatrix of  
$\vek E_{1}^{(j)}\vek S$. 
This implies that $\CH_{j-1}^{(j-1)}\cdots \CH_{1}^{(j-1)}\vek 
E_{1}^{(j-1)}\vek S$ is contained as the upper block in 
$\CH_{j-1}^{(j)}\cdots \CH_{1}^{(j)}\vek E_{1}^{(j)}\vek S$
where we recall that this sequence of reflectors was already applied at step $j-1$. Thus we
only need to apply one new reflector at iteration $j$.   
The reflector $\vek \CH^{(j)}_{j}$ only affects rows $j$ to $j+p$ of 
$\vek \CH^{(j)}_{j-1}\cdots\vek \CH^{(j)}_{1}\vek E^{(j)}_{1}\vek S$.  
This yields the
relation $\vek Z_{j} = \brac{\begin{matrix}\vek Z_{j-1} \\ \vek z_{j}^{T}
\end{matrix}}$ where $\vek z_{j}\in\R^{p}$, and we can update $\vek X_{j}$ progressively as 
an update of $\vek X_{j-1}$,
\begin{equation}\label{eqn.approx-update}
	\vek X_{j} = \vek X_{j-1} + \vek m_{j}\vek z_{j}^{T}.
\end{equation}
 
 Rather than individually storing the Householder reflector from the most recent $2p$ columns, 
 one can employ the idea presented 
in \cite{Gutknecht2008}. The authors suggested that one can store the actions
 of the Householder reflectors 
 for a block of columns as a single matrix for the 
 purpose of  applying them at future iterations. 
   This dense matrix-matrix multiplication can be performed as a level-3 
  BLAS operation.\footnote{This accumulation of the 
actions of the Householder reflections is not currently implemented in our code.}

\section{Removal of Dependent Lanczos Vectors}\label{section.dep-rhs}
We now describe some strategies for handling the linear dependence 
of a block Lanczos vector.  To maintain block size we advocate replacing
the dependent vector by a random vector, orthogonalized against all 
previous Lanczos vectors.  A set of such vectors is maintained in memory,
serving as a {\em dynamic substitutes bench} to be used upon generation of
a dependent Lanczos vector.  This procedure is described in greater detail
in Subsection \ref{subsection.dyn-sub}.

In a single-vector Krylov subspace method, 
it may happen that at step $j$, we have
that $\vek A\vek v_{j}\in\K_{j}(\vek A,\vek r_{0})$.  This implies that the grade of $\vek A$ with
respect to $\vek r_{0}$ is $j$, i.e.,
$\nu(\vek A,\vek r_{0}) = j$.
In other words, when the process creates a dependent vector the grade has been achieved.
Since $\vek r_{0}$ is the initial residual, we have that the approximation
\begin{equation}\nonumber
	\vek x_{j} = \vek x_{0} + \vek t_{j}\mbox{\ \ with\ \ }\vek t_{j}\in\CK_{j}(\vek A,\vek r_{0})
\end{equation}
is the exact solution for any Krylov subspace method (derived through a Petrov-Galerkin condition).
This situation is referred to as \textit{happy breakdown}
since it means that the true solution is contained in the existing Krylov subspace, 
see, e.g., \cite[Section 6.5.4]{Saad.Iter.Meth.Sparse.2003}.  

The notion of Krylov subspace grade has been extended to the block Krylov subspace setting \cite{GS.2009},
where we denote $\nu(\vek A,\vek F_{0})$ the block grade of $\vek A$ with respect to $\vek F_{0}$ as the
smallest integer such that 
\begin{equation}\nonumber
	\CK_{\nu(\vek A,\vek F_{0})}(\vek A,\vek F_{0}) = \CK_{\nu(\vek A,\vek F_{0})+1}(\vek A,\vek F_{0})
\end{equation}
As in the single-vector case if the block grade is achieved during the iteration of a block Krylov subspace
method, then the the method converges (if the initial block residual is used to generate the subspace).
However, unlike the single-vector case the encounter of a dependent 
basis vector in a block method does not signify that the block grade 
has been reached.
Thus if we encounter a dependent vector, this does not 
necessarily signify that the method has
converged.  The process may generate dependent 
vectors without
convergence any of the systems \cite{GS.2009}.  In the case of algorithms that are built upon the symmetric or 
nonsymmetric block/band Lanczos methods,
this dependence of the Lanczos vectors can lead to unstable algorithms if not properly handled; see, e.g., 
\cite{freund.QMR.1991,OLeary1980}.

Various strategies have been proposed to mitigate the dependence problem.  For block-level algorithms 
one must first compute or estimate the range of the block Krylov subspace basis to detect rank deficiency.
For symmetric Lanczos-based methods, O'Leary \cite{OLeary1980} 
advocates removal of the dependent vector, reducing the block size.  The update procedures for the 
systems not associated with 
the removed vector do not change, and a progressive update formula can be derived for the systems 
associated with removed
right-hand sides.  Baglama \cite{B.2000} suggests that instead of simply removing the dependent vector 
and reducing block size,
one can replace the dependent vector with a random one which has been orthogonalized against
all previous Lanczos vectors and continue unabated. For nonsymmetric Lanczos-based block QMR, 
Aliaga et al. \cite{ABFH.2000} propose to remove basis vectors before exact 
dependence is detected.  Due to 
issues of stability in
block nonsymmetric Lanczos based methods, the authors advocate defining a tolerance $d_{tol} > 0$.  
After a vector $\vek v$ has
been biorthogonalized, we have $\norm{\vek v}<d_{tol}$. We then consider $\vek v$ as almost being 
dependent, and it is removed 
from the basis. 
In \cite{ABFH.2000} a bookkeeping scheme is presented to keep track of such removals 
so that the block QMR algorithm
can be adjusted accordingly.  Recently, this technique was extended to a block conjugate gradient method 
for shifted
linear systems \cite{BF.2013}.  The bookkeeping scheme allows for the dependent vectors to be 
removed from
the process but temporarily retained in memory for the purposes of orthogonalization.

Dubrulle \cite{D.2001} proposes an alternative to the removal of dependent or near-dependent vectors for 
use in a block conjugate gradient algorithm.  He
proposes to use a change-of-basis strategy for the block descent directions and other algorithmic changes 
to avoid the problem
long before near-rank deficiency of the block basis vectors occurs.  This additionally avoids the need for 
basis rank estimation.

In \cite{RS.2006} following from \cite{L.2003} the authors suggest that removing nearly dependent
directions could represent an unacceptable loss of information.  They recommend instead to
reintroduce the dependent directions at the next iteration.  They also consider some 
different methods of defining and detecting near breakdown.

For our version of the block MINRES algorithm we define dependence as in  \cite{ABFH.2000, BF.2013} 
with the candidate vector being considered dependent if $h_{j+p,j} \leq \gamma$ for an  a priori chosen
$0 \leq\gamma\ll 1$.  One of the characteristics of a block Krylov method built from the banded 
Lanczos method is that 
there is no need for any basis rank estimation.  Since we construct only one band Lanczos vector at 
a time, we simply need
to compute $h_{j+p,j}$, i.e., compute the norm of the newest basis vector after orthogonalization via 
the band Lanczos 
process.  This has been observed previously (in the context of eigenvalue computations) \cite{B.2000}.  
Baglama presents two options for dealing with linear dependence.  One option is to reduce 
block size by one and adjust short-term recurrences accordingly.  The other is to generate a random vector and
orthogonalize it with respect to all previous Lanczos vectors.  This normalized vector is then put in the place of
the dependent Lanczos vector.  Either option results in minimal changes to  
the algorithm, but
we only discuss the algorithmic modifications required to incorporate the latter,
as we favor maintaining block size.
For illustration, we present examples for a particular block size, for ease of discussion; 
however, it is clear
that the simplifications presented do not change if the block size increases.

To see a discussion of the algorithmic ramifications of adopting the block size 
reduction strategy, see the research report \cite{S.2014-arxiv}.
If we assume for simplicity that we only remove truly dependent vectors, i.e., 
$\gamma=0$.
In \cite{S.2014-arxiv} it is shown that for 
each dependent vector removed we will have a two-vector reduction in storage requirements 
for the construction
of the search directions.  In total, for each block size reduction, we have a four-vector reduction in 
storage requirements.

\subsection{Maintaining block size using a dynamic substitutes bench}\label{subsection.dyn-sub}
We now discuss how inserting a random orthogonalized vector into the basis
affects the algorithm.
We then present a strategy for having random orthogonal vectors available.

We begin by describing the replacement procedure in more detail.  At iteration $j$, we compute 
$\vek A\vek u_{j}$.  After orthogonalization
we see that $h_{j+p,j} = 0$.  Thus, $\vek A\vek u_{j}$ is in the range of the previous 
Lanczos vectors.
Let $\widehat{\vek u}_{j+p}$ be a vector constructed by taking a random vector $\widehat{\vek w}$, 
orthogonalizing $\widehat{\vek w}$
with respect to all previous Lanczos vectors and setting  
$\widehat{\vek u}_{j+p} = \widehat{\vek w}/\norm{\widehat{\vek w}}$.
Then the algorithm continues as before with this modified block Lanczos basis.  

This strategy allows us to maintain the block size $p$ when a loss of independence is encountered.
We advocate this policy specifically in the context of high-performance computing applications.
Of course, this must be weighted against the costs of maintaining a larger block size.
If we generate the block Krylov subspace for $p$ initial residuals to solve (\ref{eqn.AXB}), 
$\K_{k,m}(\vek A,\vek F_{0})$ will contain sufficient information to construct high-quality solutions for all 
right-hand sides for large enough values of $j=(k-1)p+m$ in theory \cite{GS.2009}.
However, in practice exact convergence in this scenario would not occur.  Maintaining the large
block size allow us to build a larger constraint space for each block matrix-vector product executed.  
In the high-performance computing setting the low costs of this strategy make them 
worthy of consideration.

What modifications must be made to the block MINRES algorithm to accommodate this strategy?  
It turns out, very few.  Of course, 
we do not store the complete Lanczos basis, as this would defeat the purpose of developing a method 
for symmetric systems.
However, we need to orthogonalize the random vector against the entire basis.  
As a work-around we can generate a random vector at the start of the iteration and simply orthogonalize 
against each Lanczos vector 
as it is created.  This would require only one additional vector of storage and an additional orthogonalization 
per iteration.
If we are solving a problem in which we expect there to be more than one occurrence of loss of 
linear independence then
we can generate more than one random vector, balancing between increasing the storage requirements 
and insuring against
the dependence problem.  This strategy does entail additional computational cost, 
but it allows us to achieve the goal of maintaining block size $p$ when breakdown occurs,
and it is desirable to maintain the larger block size for the data movement efficiencies previous 
discussed.

One might be concerned that introducing a vector not created by the band Lanczos process will 
destroy the short-term recurrences
which make symmetric Lanczos methods so attractive.  However, this is not the case. Suppose 
that after iteration $j$ 
we continue the band Lanczos process
with the modified basis.  Let $\widehat{\vek U}_{j+p}\in\R^{n\times (j+p)}$
be the matrix containing the band Lanczos vectors but with $\widehat{\vek u}_{j+p}$ as its last column.  
Observe that the matrix $\widehat{\vek U}_{j+p}^{\ast}\vek A\widehat{\vek U}_{j+p}$ is symmetric; 
inserting the new basis vector
does not affect this.  Thus the banded structure of $\widehat{\vek H}_{j}$ defined by 
$\vek A\widehat{\vek U}_{j}=\widehat{\vek U}_{j+p}\widehat{\vek H}_{j}$ is 
the same as that of $\overline{\vek H}_{j}$.  The only change is that we now have zero entries at 
$h_{j+p,j}$ and $h_{j,j+p}$.
This in turn gives a slight change in structure to $\widehat{\vek R}_{j}$, the upper triangular factor 
in the QR-factorization of $\widehat{\vek H}_{j}$.

As an example, suppose $p=2$ and that $\vek A\vek v_{5}$ is in the span of the existing band Lanczos 
vectors, as in the last example.  If  
we continue the band Lanczos process with the modified basis,
we have the following structures for $\widehat{\vek H}_{8}$ and 
\linebreak $\widehat{\vek R}_{8}\in\R^{10\times 8}$ (right), 
\begin{equation}\nonumber
\resizebox{.49\hsize}{!}{$\begin{bmatrix}  h_{1,1}&h_{1,2}&h_{1,3}&	       &  		&	      &			&	         &	        &\\
				               			   h_{2,1}&h_{2,2}&h_{2,3}&h_{2,4}&	   	&   	      &			&	         &	        &\\
				   				   h_{3,1}&h_{3,2}&h_{3,3}&h_{3,4}&h_{3,5}	&	      &			&	         &	        &\\
				 					       &h_{4,2}&h_{4,3}&h_{4,4}&h_{4,5}	&h_{4,6}&			&	         &	        &\\
				 					       &	    &h_{5,3}&h_{5,4}&h_{5,5}	&h_{5,6}&\mathbf{0} &               &	        &\\
				 					       & 	    &		 &h_{6,4}&h_{6,5}	&h_{6,6}&h_{6,7}	&h_{6,8}  &	        &\\
				 					       &	    &		 &	       &\mathbf{0}&h_{7,6}&h_{7,7}	&h_{7,8}  &h_{7,9}   &\\
				  					       &	    &		 &	       &		&h_{8,6}&h_{8,7}	&h_{8,8}  &h_{8,9}   &h_{8,10}\\
				     					       &  	    &		 &	       &		&	      &h_{9,7}	&h_{9,8}  &h_{9,9}   &h_{9,10}\\
				  					       &	    &		 &	       &		&	      &		   	&h_{10,8}&h_{10,9}&h_{10,10}\\
				  					       &	    &		 &	       &		&	      &		   	&               &h_{11,9}&h_{11,10}\\
				  					       &	    &		 &	       &		&	      &		   	&               &               &h_{12,10}
				            \end{bmatrix}$}\mbox{\ \ }
				            \resizebox{.49\hsize}{!}{$\begin{bmatrix}  r_{1,1}&r_{1,2}   &r_{1,3} &r_{1,4}&r_{1,5}	&	     &			&	     &		         &\\
				               			  	       &r_{2,2} &r_{2,3} &r_{2,4}&r_{2,5}	&r_{2,6}&			&            &		         &\\
				   				   	       &	    &r_{3,3} &r_{3,4}&r_{3,5}	&r_{3,6}&\mathbf{0}  &            &		         &\\
				 					       &	    &		 &r_{4,4}&r_{4,5}	&r_{4,6}&r_{4,7}	&r_{4,8}&		         &\\
				 					       &	    &		 &	       &r_{5,5}	&r_{5,6}&r_{5,7}	&r_{5,8}&\mathbf{0} &\\
				 					       & 	    &		 &	       &		&r_{6,6}&r_{6,7}	&r_{6,8}&r_{6,9}       &r_{6,10}\\
				 					       &	    &		 &	       &		&	      &r_{7,7}	&r_{7,8}&r_{7,9}       &r_{7,10}\\
				  					       &	    &		 &	       &		&	      &		   	&r_{8,8}&r_{8,9}       &r_{8,10}\\
				     					       &  	    &		 &	       &		&	      &		   	& 	     &r_{9,9}       &r_{9,10}\\
				  					       &	    &		 &	       &		&	      &		   	&	     &		        &r_{10,10}\\
				  					       &	    &		 &	       &		&	      &		   	&	     &		        &\\
				  					       &	    &		 &	       &		&	      &		   	&	     &		        &
				            \end{bmatrix}$}.
\end{equation}
Observe that due to the symmetry of $\widehat{\vek H}_{8}$ (left), 
$\norm{\vek v_{6}}=0$ gives us an additional
zero in the super diagonal. This yields two zero entries in 
the upper-most superdiagonal of 
$\widehat{\vek R}_{8}\in\R^{10\times 8}$ (right).
This indicates that the final effects of replacing the dependent vector with a random one are 
minimal.  The two zeros
are introduced into upper Hessenberg matrix but the bandwidth and symmetry properties remain unchanged.  The
introduction of a zero in the seventh column of $\widehat{\vek R}_{8}$ and another in the ninth simply means that 
the seventh and ninth band Lanczos vectors are linear combinations of the previous four rather than the 
previous five
search directions, recalling the construction of the search directions (\ref{eqn.search-dir-comp}).

In this discussion we have assumed exact 
deflation, i.e., the newest candidate Lanczos vector is exactly in the span of the previous vectors.  In 
practice we want to remove a generated vector when it is "nearly" dependent, i.e., reject
$\vek v_{j}$ when $h_{j+p,j} < \gamma$ where $\gamma<1$ is some dependence tolerance constant
sufficiently far from zero, 
as in \cite{BF.2013,FM.1997}.
This is especially true in block Krylov method formulations relying on short-term recurrences that
are formulated with a progressive update of the solution at each iteration.
In our code, any removed vector is held in storage and new Lanczos vectors are still orthogonalized 
with respect to it until the removed vector would naturally have been dropped due to the band Lanczos
relation.
It should be noted; removing basis vectors
in this way no longer follows the mathematical derivation, and we must understand the effect of this strategy
on convergence, choosing $\gamma$ in a way that balances our need for stability with any delay in 
convergence this strategy might cause.

\section{Convergence Theory}\label{section.conv}
Theoretically block MINRES is a version of block GMRES for symmetric systems.  
Simoncini and Gallopoulos discussed the convergence properties of block 
GMRES \cite{Simoncini.Conv-Block-GMRES} including a result by Vital \cite{Vital1990}.  
We can describe the quality of the residual produced at iteration $j$.  
For $\vek b^{(i)}$ the $i$th column of the right-hand side $\vek B$, Algorithm 
\ref{alg.block-minres} minimizes the $i$th column of the residual $\vek f_{j}^{(i)}$ over the 
subspace $\K_{k,m}(\vek A,\vek F_{0})$.  Thus we can expect $\norm {\vek 
f_{j}^{(i)}}$ to be at least as good as the norm of the residual produced by running $J$ 
steps of MINRES with ${\vek b^{(i)}}$ as the single right-hand side where 
\begin{equation}\nonumber
	J = \begin{cases}k+1\text{ if } i \leq m\\ k\text{ if } i > m\end{cases}.
\end{equation}  
This easily can be 
understood by recalling the definition of $\K_{k,m}(\vek A,\vek F_{0})$ in 
(\ref{eqn.block-Krylov-subspace}).  We observe that having a larger subspace over which to 
minimize is not guaranteed to give improvements in convergence. The additional information contained in $
\K_{k,m}(\vek A,\vek F_{0})$ may not be helpful in the minimization process.  
For specially related 
right-hand sides, we may have a great boost in performance.

In theory, a block Krylov subspace iterative method may terminate before the subspace becomes the
full space $\Rn$.  As in the single vector case, achievement of the block grade implies that the exact 
solution correction is in that subspace. 
For a block Arnoldi (or Lanczos)-based method, if 
$\nu(\vek A,\vek F_{0})$ is the block grade then 
\begin{equation}\nonumber
	\vek X = \vek X_{0} + \vek T_{\nu(\vek A,\vek F_{0})} \mbox{\ \ for\ \ } \vek T_{\nu(\vek A,\vek F_{0})}\in \K_{\nu(\vek A,\vek F_{0})}(\vek A,\vek F_{0})
\end{equation}
is the exact solution, where the correction is generated by any Petrov-Galerkin condition
(since $\vek F_{0}$ is the initial residual).
Using the equivalence between the block Arnoldi and band Arnoldi bases at
specific iterations (\ref{eqn.block-band-equiv}), we can extend this notion of block grade
to a block Krylov subspace generated by a band rather than truly block process.
From (\ref{eqn.block-band-equiv}) and the definition of block grade, 
we have the following straightforward identity,
\begin{proposition}\label{prop.block-grade-ident}
	Let $\nu(\vek A,\vek F_{0})$ be the block grade of $\vek A$ and $\vek F_{0}$, as defined in 
	{\rm\cite{GS.2009}}.  Then we have 
	\begin{equation}\nonumber
		\K_{\nu(\vek A,\vek F_{0})}(\vek A,\vek F_{0}) = \K_{\nu(\vek A,\vek F_{0}),0}(\vek A,\vek F_{0})
	\end{equation}
	where the band Krylov subspace satisfies the invariance
	\begin{equation}\label{eqn.banded-invariance}
		 \K_{\nu(\vek A,\vek F_{0}),0}(\vek A,\vek F_{0}) =  \K_{\nu(\vek A,\vek F_{0})+1,0}(\vek A,\vek F_{0}).
	\end{equation}
\end{proposition}

We must be careful when describing the notion of grade for a band Krylov method.  
Depending on the ordering of the right-hand sides, we will have that at some iteration $K$ there 
will be no further increase of the band Krylov subspace dimension where
\begin{equation}\nonumber
	K\in\curl{(\nu(\vek A,\vek F_{0})-1)p+1,(\nu(\vek A,\vek F_{0})-1)p+2,\ldots, (\nu(\vek A,\vek F_{0}))p}.
\end{equation}
However, the exact value of $K$ depends on the ordering of the right-hand sides, since each has 
an associate grade with respect to the matrix $\vek A$, and from \cite[Lemma 6]{GS.2009} these 
individual grades can be related
to the block grade by  
\begin{equation}\nonumber
	 \K_{\nu(\vek A,\vek F_{0}),0}(\vek A,\vek F_{0}) = \sum_{i=1}^{p}\CK_{\nu_{i}}(\vek A,\vek f_{0}^{(i)})\mbox{\ \ where\ \ }\nu_{i} = \nu(\vek A,\vek f_{0}^{(i)}).
\end{equation}
Thus to describe the
notion of grade in an unambiguous way, we must make one assumption.  Without loss of generality 
we assume that the columns of $\vek F_{0}$ are ordered such that the single-vector Krylov
subspace grades with respect to each right-hand side satisfy
\begin{equation}\label{eqn.rhs-grade-ordering}
	\nu(\vek A,\vek f_{0}^{(1)})\leq \nu(\vek A,\vek f_{0}^{(2)}) \leq \cdots\leq \nu(\vek A,\vek f_{0}^{(p)}).
\end{equation}
This determines the order in which the individual grades are achieved 
in (\ref{eqn.block-Krylov-subspace}) as we iterate.
By fixing the ordering of the columns of $\vek F_{0}$ as in 
 (\ref{eqn.rhs-grade-ordering}), the iteration $K$ at which we have achieved
 the largest possible block Krylov subspace dimension (constructed by a band Arnoldi-based
  method) can be unambiguously defined as the \textit{band grade} with respect to $\vek A$ and $\vek F_{0}$.  
  Thus the notion
of block grade described in \cite{GS.2009} can be translated unambiguously to a
band-Arnoldi based method.

\section{A Note on Implementation}\label{section.impl}

We conclude our description with some notes about practical implementation details.  
In order to achieve the data movement benefits of block operations we apply the operator to a 
block of $p$ vectors every $p$ iterations.
We must store $2p$ Lanczos vectors, 
$2p$ search directions, $2p$ Householder reflections, the lower subdiagonal entries of $p$ previous columns of $\overline{\vek H}_{j}$, and the $j$th column of 
$\overline{\vek H}_{j}$.    We also may store some nearly dependent vectors for orthogonalization
and some random vectors used to replace dependent vectors.

While the symmetry of $\vek A$ allows for a fixed storage requirement
we must take care with how we store the Lanczos vectors and search directions.
Our primary goal in describing storage layout is to show how the method can be implemented 
without much need for tracking of indices.
For simplicity of implementation we advocate that the Lanczos vectors and search directions be 
stored in a first-in-first-out (FIFO) queue holding $2p$ vectors.  
This results in the most recently generated vector will be in the last position in the queue; and 
when a new vector is created, the oldest vector will automatically be overwritten.

The full matrix $\overline{\vek H}_{j}$ need not be stored, 
but the lower subdiagonal entries are needed for the block
Lanczos process (as they are orthogonalization coefficients in future iterations due to symmetry).
The subdiagonal entries from the $p$ most recent columns of $\overline{\vek H}_{j}$ can be stored
in a FIFO queue (or in a $p\times p$ matrix called $\vek C_{p\times p}$ behaving as a queue
with the newest entries inserted into the last column).  
Storing the entries in this manner results in the nonzero superdiagonal entries of the current column of 
$\overline{\vek H}_{j}$ being available as the 
nonzero antidiagonal entries of $\vek C_{p\times p}$.  This allows us to obtain the super diagonal entries of the 
current column without computing the associated inner 
products.  

For block size $p=5$ at iteration $j=7$ of the banded
Lanczos process, we have
\begin{equation}\label{eqn.anti-diag-symmetry}
	\overline{\vek H}_{7}(:,7) = \brac{\begin{matrix}0 \\ {\textbf h_{2,7}} \\ {\textbf h_{3,7}} \\ {\textbf h_{4,7}} \\ {\textbf h_{5,7}} \\ {\textbf h_{6,7}} \\ h_{7,7}\\ h_{8,7} \\ 
h_{9,7} \\ h_{10,7} \\ h_{11,7} \\ h_{12,7} \\ h_{13,7} \\ h_{14,7}\end{matrix}}\mbox{\ \ \ and\ \ \ } \vek C_{p\times p} = \brac{\begin{matrix}h_{3,2}&h_{4,3}&h_{5,4}&h_{6,5}
&{\textbf h_{7,6}}\\ h_{4,2}&h_{5,3}&h_{6,4}&{\textbf h_{7,5}}&h_{8,6}\\ h_{5,2}&h_{6,3}&{\textbf h_{7,4}}&h_{8,5}&h_{9,6}\\ h_{6,2}&{\textbf h_{7,3}}
&h_{8,4}&h_{9,5}&h_{10,6}\\ {\textbf h_{7,2}}&h_{8,3}&h_{9,4}&h_{10,5}&h_{11,6}\end{matrix}}.
\end{equation}
In (\ref{eqn.anti-diag-symmetry}), note the correspondence between bold entries in 
$\overline{\vek H}_{7}(:,7)$ and 
the antidiagonal entries of $\vek C_{p\times p}$, computed at previous iterations.

\begin{algorithm}[b!]
\caption{Block MINRES (band Lanczos Version) with no Breakdown}
\label{alg.block-minres}

\SetKwInOut{Input}{Input}\SetKwInOut{Output}{Output}
\Input{$\vek A\in \Rnn$ Symmetric, $\vek B\in\R^{n\times p}$, $\vek X_{0}=\vek 0$ , $\epsilon > 0$, $M\in\N$ }
\Output{$\vek X\in\R^{n\times p}$ such that $\norm{\vek B(:,j) - \vek A\vek X(:,j)}/\norm{\vek B(:,j) - \vek A\vek X_{0}(:,j)}<\epsilon\ \forall j \leq p$}

Compute the QR-Factorization $ \vek B \vek V_{p}\vek S$\\
$\widehat{\vek S} \leftarrow \vek S\vek E_{1}^{(1)}$\\
$\vek X \leftarrow \vek X_{0}$\\
$\vek R \leftarrow \vek B - \vek A\vek X$\\
\While{$\max_{0<i\leq p}\curl{\norm{ (\overline{\vek Z}_{j}^{(i)})_{j+1:j+i} }} < \varepsilon\norm{\vek b^{(i)}}$ \mbox{\rm and} $j \leq M$}{
	\If{$j\mod p = 1$}{
		$\vek W \leftarrow \vek A\begin{bmatrix} \vek v_{j} & \vek v_{j+1} & \cdots & \vek v_{j+p}  \end{bmatrix}$\\
	}
	\If{$j\mod p \neq 0$}{
		$\ell = j\mod p$\\
	}\Else{
		$\ell = p$\\
	}
	$\vek w \leftarrow \vek W(:,\ell)$\\
	\If{$j > 1$}{
		\For{$i=j-p:j-1$}{
			$h_{i,j} = h_{j,i}$\\
			$\vek w \leftarrow h_{i,j}\vek w$\\
		}
	}
	\For{$i=j:j+p-1$}{
		$h_{i,j} = \vek v_{i}^{\ast}\vek w$\\
		$\vek w \leftarrow h_{i,j}\vek v_{i}$\\
	}
	$h_{j+p,j} = \norm{\vek w}$\\
	$\vek v_{j+1} = \vek w/h_{j+p,j}$\\
	\If{$j>1$}{
		$\overline{\vek r}_{j}^{(j)} \leftarrow \CH_{j-1}^{(j)}\cdots\CH_{j-2p}^{(j)}\overline{\vek h}_{j}^{(j)}$\\
	}
	Generate Householder reflection $j$th column of $\overline{\vek H}_{j}$\\
	$\overline{\vek r}_{j}^{(j)} \leftarrow \CH_{j}^{(j)}\overline{\vek r}_{j}^{(j)}$\\
	\If{m = 1}{
		$\vek m_{1} = \vek v_{1}/\overline{\vek R}_{1}(1,1)$\\	
	}\Else{
		$\vek w \leftarrow \vek v_{j}$\\
		\For{$i = j-2p:j-1$}{
			$\vek w \leftarrow \vek w - \overline{\vek R}_{j}(i,j)\vek m_{i}$\\
		}
		$\vek m_{j} = \vek w/\overline{\vek R}_{j}(j,j)$\\
	}
	$\vek z^{T} \leftarrow \widehat{\vek S}(j,:)$\\
	$\vek X \leftarrow \vek X + \vek m_{j}\vek z^{T}$\\
	$\widehat{\vek S} \leftarrow \brac{\begin{matrix}\widehat{\vek S} \\ \vek 0_{1\times p}\end{matrix}}$\\
	$j\leftarrow j+1$\\
}
\end{algorithm}

\section{Numerical Results}\label{section.num}

We present numerical experiments to demonstrate the effectiveness and behavior of Algorithm 
\ref{alg.block-minres}.   In 
all experiments, we compared the performance of block MINRES with sequential applications of Matlab's 
MINRES function.  We compared performance 
using iteration counts and sometimes CPU timings.    
However, note that if we measure the cost of an iteration according to a data movement
metric, the cost of the iteration would be dominated by the block matrix-vector product executed every $p$ 
iterations, amortized over the subsequent $p$ iterations.  The block matrix-vector product does not
cost (in data movement) $p$ times as much as $p$ single matrix-vector products \cite{PSS.2013}.
In this metric, an iteration of our method and a sequential MINRES iteration are not equivalent.   

All tests were performed on a Macbook Pro containing a 2.3 GHz Intel 
Core i5 processor with 8 GB of 1333MHz DDR3 main memory running the 64-bit version of Matlab R2011b. In 
any experiment involving the generation of random vectors, we used Matlab's \verb|mt19937ar| 
random number generator,
with seed $0$, which was initialized at the beginning of each experiment.  
The tests were performed for a model shifted Laplacian problem.
Let $\vek L\in\R^{n_{1}\times n_{1}}$, with $n_{1} = 40000$, 
be the discretization of the Laplacian operator on a $200\times 200$ regular grid using central differences,
constructed by setting $\vek T = \verb|tridiag(1,-2,1)|$ and $\vek L=h^{-2}(\vek I\oplus\vek T + \vek T\oplus\vek I)$
where $h=1/199$.  This 
matrix is negative-definite.    Let $\vek A = -\vek L - 200\vek I$.  Due to the eigenvalue distribution of $\vek L$, 
we have that $\vek A$ is indefinite.  In all experiments, we precondition with the incomplete Cholesky
factors of $-\vek L$ constructed using Matlab's \verb|ichol()| function with the default settings.  

\begin{figure}[htb]
\hfill
\includegraphics[scale=0.40]{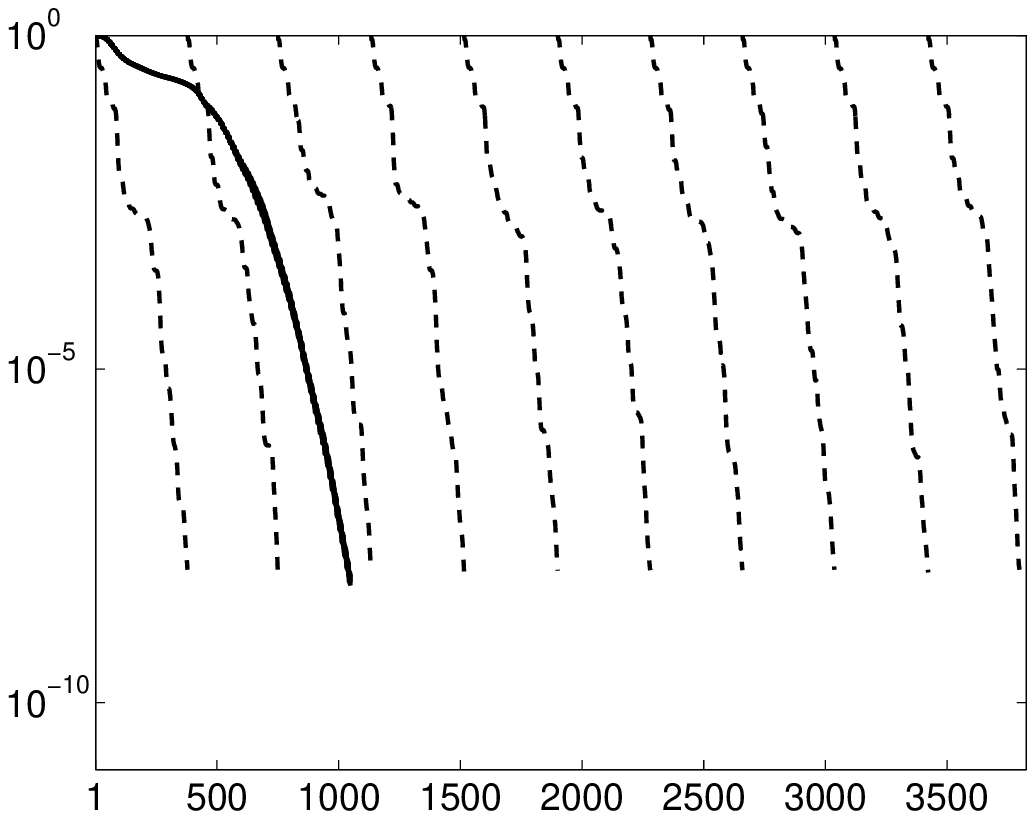}\quad\quad\hfill
\vspace*{2em}
\begin{picture}(0,0)
\put(-268,-15){{\footnotesize Block MINRES: 1048 iter., 16.798 sec.}}
\put(-268,-25){{\footnotesize Sequential MINRES: 3787 iter., 23.546 sec.}}
\end{picture}
\caption{Comparison of the performance of Algorithm {\rm \ref{alg.block-minres}} versus sequential applications of
MINRES on the discretized Laplacian system with ten randomly generated right-hand sides.  
The {\em solid black} curve
is actually the ten convergence curves for each right-hand side when solved by Algorithm {\rm\ref{alg.block-minres}}
overlaid on one another.  We see that
in the case of these ten right-hand sides that block MINRES convergence for all ten systems is 
qualitatively the same.  The
{\em black dashed} curves are the convergence curves for each sequential application of MINRES 
for each right-hand side.}
\label{fig.blockMINRES-10rhs}
\end{figure}

We begin by demonstrating the performance of the algorithm on the shifted Laplacian system 
with ten randomly generate right-hand sides.  
In Figure \ref{fig.blockMINRES-10rhs}, we see that for these right-hand sides, the block MINRES algorithm 
converges in fewer iterations and less time.

\begin{figure}[htb]
\hfill
\includegraphics[scale=0.40]{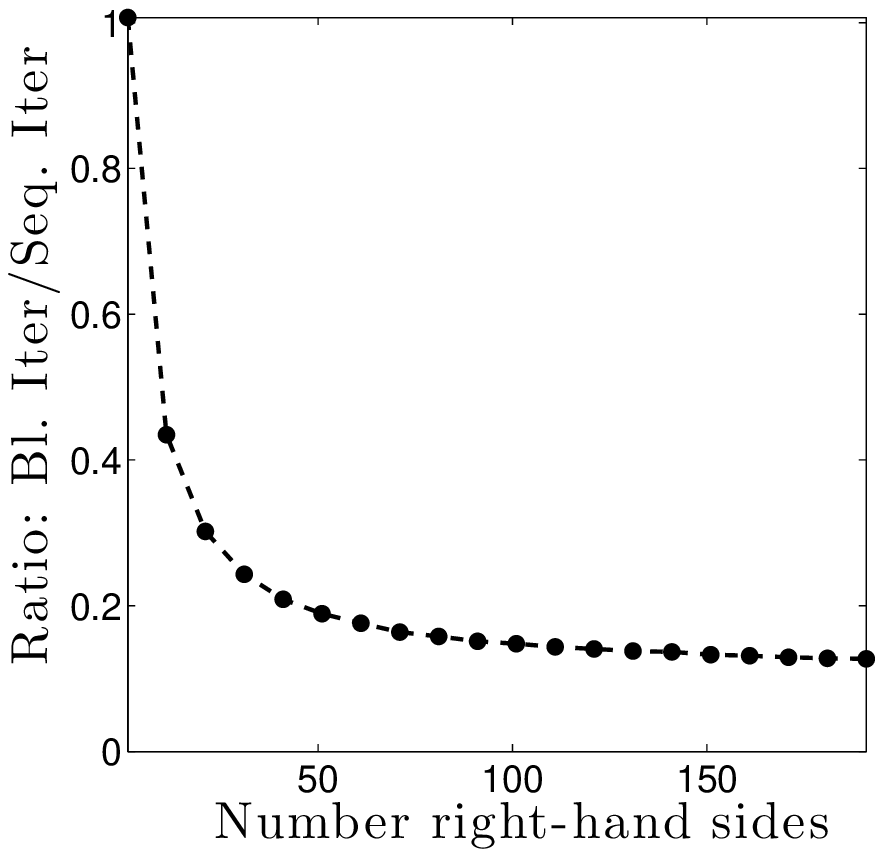}\quad\quad\hfill
\vspace*{-1em}
\caption{For different numbers of right-hand sides, we plot the ratio between the number
of iterations required by our block MINRES method and the sequential MINRES method to solve the
system with multiple right-hand sides to a tolerance of $10^{-8}$.  As the number of right-hand
sides increases, the ratio decreases, i.e., our method requires fewer iterations than sequential
MINRES.}
\label{fig.blockMINRES-rhsNumTest-iter}
\end{figure}

We can also compare performance of our method versus sequential applications of MINRES
for varying numbers of right-hand sides.  We take as our first right-hand side the vector of 
all ones.  If we have $p$ total right-hand sides, we take the remaining $p-1$ to be the first
$p-1$ columns of the $\vek I_{n_{1}}$.  In Figure \ref{fig.blockMINRES-rhsNumTest-iter}, 
we plot for various $p$, the ratio between
the iteration count of our method and the total iteration count for $p$ sequential applications of 
Matlab's MINRES.  For this experiment, we see a reduction in the ratio as $p$ increase,
but the marginal benefit of adding each additional right-hand side diminishes for
larger numbers of right-hand sides.

\begin{figure}[htb]
\hfill
\includegraphics[scale=0.40]{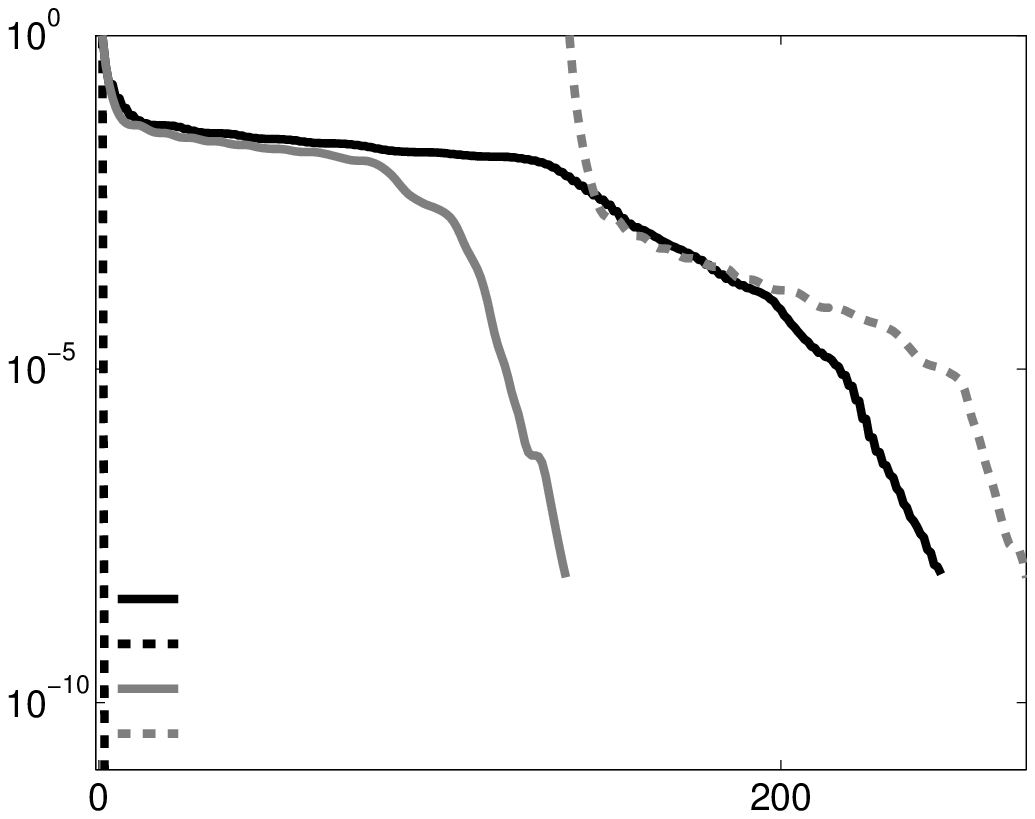}\quad\quad\hfill
\vspace*{-1	em}
\begin{picture}(0,0)
\put(-260,39){{\footnotesize Bl. MR $\vek e_{1}$}}
\put(-260,31){{\footnotesize Bl. MR $\vek A\vek e_{1}$}}
\put(-260,23){{\footnotesize Seq. MR $\vek e_{1}$}}
\put(-260,13){{\footnotesize Seq. MR $\vek A\vek e_{1}$}}
\end{picture}
\caption{Demonstration of the algorithm's performance in the case that it encounters 
dependence.  In this case, with the right-hand sides
$\vek e_{1}$ and $\vek A\vek e_{1}$, dependence occurs at the first iteration. Since the first right-hand 
side is the solution to the second system, 
we get immediate convergence for the second system, and block MINRES continues for the other system, 
replacing the dependent 
vector with a random one.}
\label{fig.blockMINRES-dep-test}
\end{figure}

We demonstrate that our removal of dependent vectors works as described.  Of course, it is 
difficult to choose a pair
of right-hand sides for which dependence will occur in later iterations.  Thus, as a simple, 
easy-to-construct test, we chose the first 
right-hand side $\vek e_{1}$, as the first canonical basis vector.  The second right-hand side is 
$\vek A\vek e_{1}$, the image of
the first canonical basis vector, i.e., the first column of our coefficient matrix.  This will result in 
dependence at the first iteration of our algorithm.
As is shown in Figure \ref{fig.blockMINRES-dep-test}, this leads to immediate convergence for 
that system when running block MINRES.  Of course, this example is 
not likely to occur in practice.  It merely demonstrates that the algorithm can handle dependence gracefully.

We demonstrate how the relationship between the right-hand sides can affect the 
performance of  
block MINRES. 
\begin{figure}[htb]
\hfill
\includegraphics[scale=0.35]{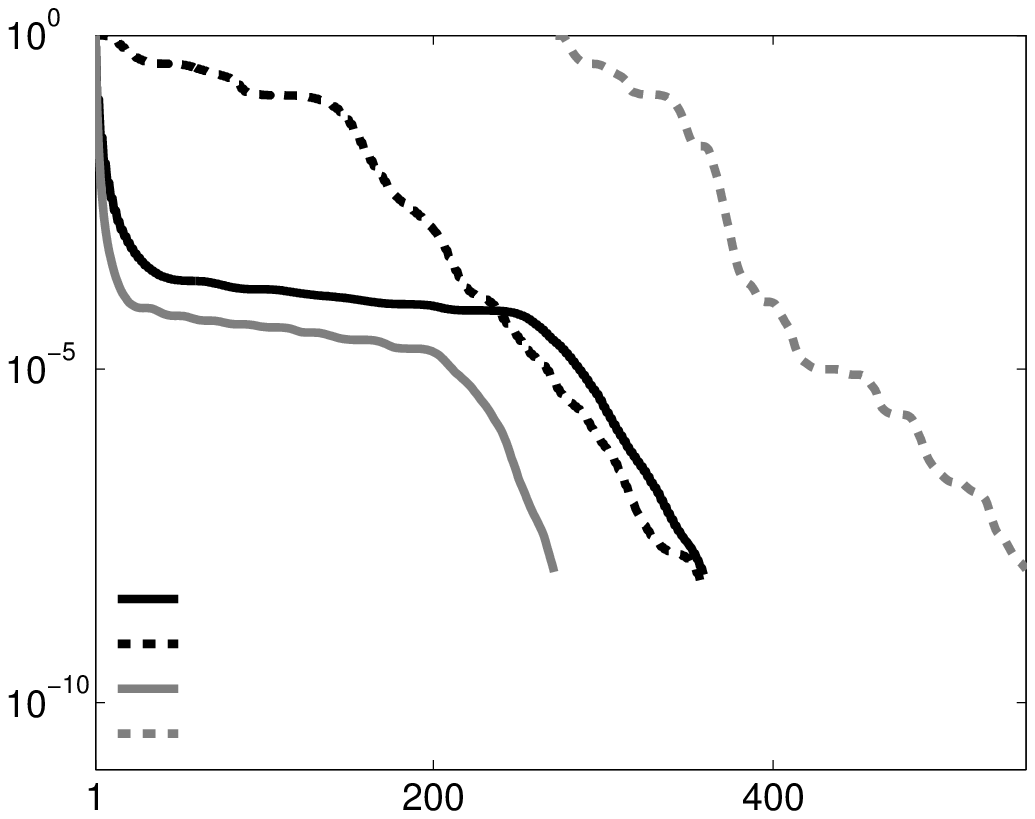}\quad\quad
\includegraphics[scale=0.35]{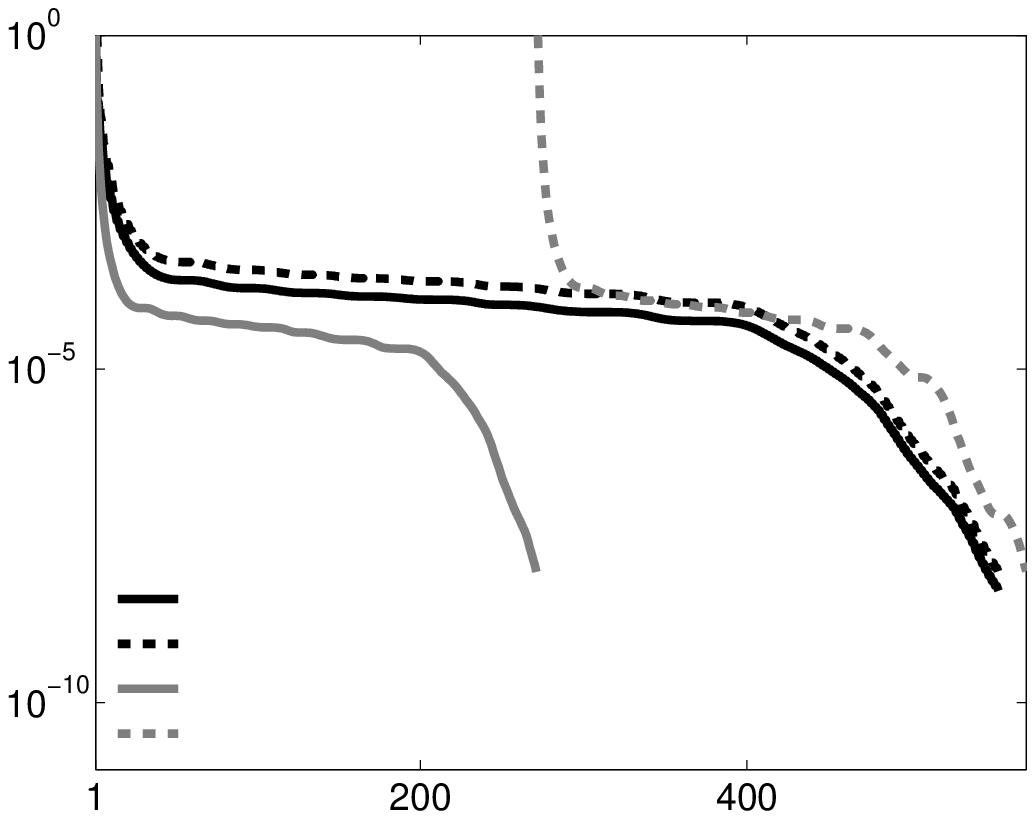}\hfill
\vspace*{1em}
\begin{picture}(0,0)
\put(32,57){{\footnotesize Bl. MR $\vek b_{1}$}}
\put(32,49){{\footnotesize Bl. MR $\vek b_{2}$}}
\put(32,41){{\footnotesize Seq. MR $\vek b_{1}$}}
\put(32,31){{\footnotesize Seq. MR $\vek b_{2}$}}
\put(20,10){{\footnotesize Block MINRES: 358 iter., 2.0934 sec.}}
\put(20,0){{\footnotesize Sequential MINRES: 547 iter., 3.1552 sec.}}
\put(230,57){{\footnotesize Bl. MR $\vek b_{1}$}}
\put(230,49){{\footnotesize Bl. MR $\hat{\vek b}_{2}$}}
\put(230,40){{\footnotesize Seq. MR $\vek b_{1}$}}
\put(230,31){{\footnotesize Seq. MR $\hat{\vek b}_{2}$}}
\put(220,10){{\footnotesize Block MINRES: 553 iter., 3.2999 sec.}}
\put(220,0){{\footnotesize Sequential MINRES: 569 iter., 3.3379 sec.}}
\end{picture}
\caption{Performance of block MINRES for different right-hand sides.  In the figure on the left, the two right-hand
 sides are  $\vek b_{1} = \vek e_{n_{1}}^{(1)}$ and $\vek b_{2} = \vek 1$.  In the figure on the right, 
 $\vek b_{1}$ does not change,
 but $\hat{\vek b}_{2} = \vek e_{n_{1}}^{(2)}$.}
\label{fig.blockMINRES-performance-diffRHS}
\end{figure}
\begin{figure}[htb]
\hfill
\includegraphics[scale=0.40]{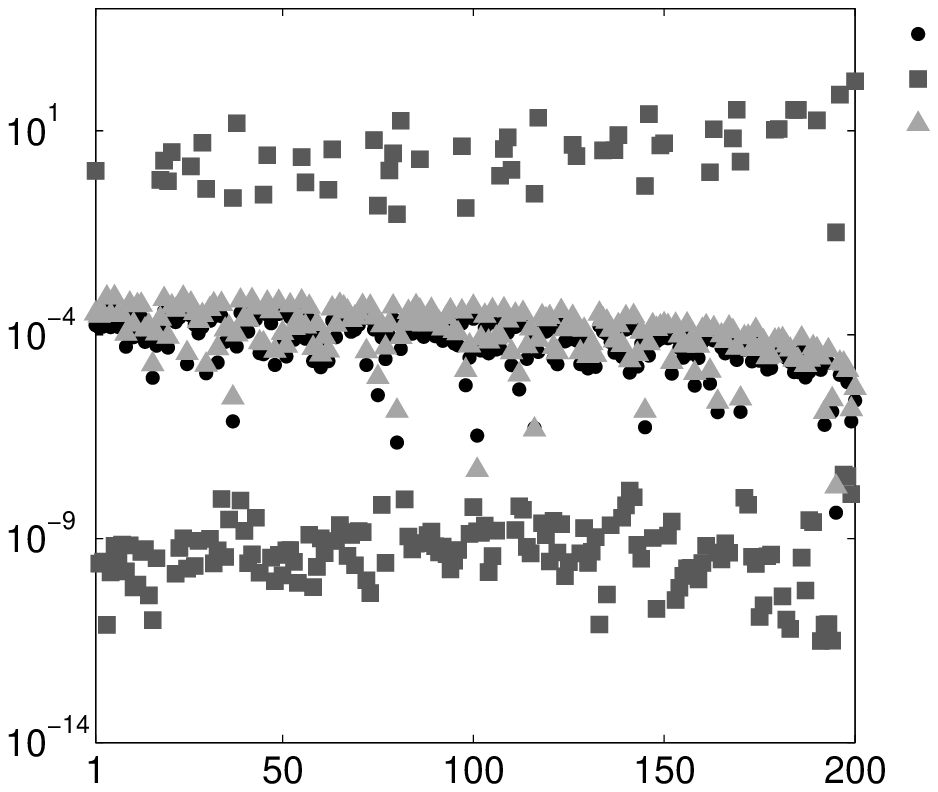}\quad\quad\hfill
\vspace*{-.9em}
\begin{picture}(0,0)
\put(-110,144){{\footnotesize $\vek b_{1}$}}
\put(-110,135){{\footnotesize $\vek b_{2}$}}
\put(-110,123){{\footnotesize $\hat{\vek b}_{2}$}}
\end{picture}
\caption{Magnitude of the components of different right-hand-sides in the eigenspace spanned 
by the two hundred eigenvectors associated with the smallest magnitude 
eigenvalues.}
\label{fig.blockMINRES-eigFig}
\end{figure}
We 
compared the performance of our block MINRES implementation with that of sequential runs of Matlab's 
MINRES for $\vek A$ with three pairs of right-hand sides.  
For the first pair, let $\vek b_{1} = \vek e_{n_{1}}^{(1)}$ and $\vek b_{2} = \vek 1$, the vector of all ones.   
For second pair of right-hand sides, we let $\hat{\vek b}_{1} 
= \vek b_{1}$ but change the second right-hand side by letting $\hat{\vek b}_{2} = \vek e_{n_{1}}^{(2)}$.  
In Figure \ref{fig.blockMINRES-performance-diffRHS}, we 
show a comparison of convergence curves for these pairs of right-hand sides.  We observe that 
exchanging $\vek b_{2}$ for $\hat{\vek b}_{2}$ degrades the 
performance of our Block MINRES implementation.  Recall that the convergence of a Krylov subspace 
method for a symmetric 
system is completely determined by its eigenvalues
and the decomposition of the initial residual in the eigenbasis.  
For an indefinite system, the eigenvalues closest to the origin cause a delay in convergence.  
In Figure \ref{fig.blockMINRES-eigFig}, we decomposed the three right-hand sides in the
eigenbasis and plotted the magnitudes of the $200$ eigencomponents associated to small 
eigenvalues.  What we see is that almost all the components of $\vek b_{1}$ and 
$\vek b_{2}$ have similar magnitude
while those of $\hat{\vek b}_{2}$ differ, with some being larger and others being smaller. 
Therefore, we 
hypothesize that a pair of right-hand sides that have strong components from different parts of the 
eigenspace might 
complement each other well.  

We concoct some experiments to explore this line of thinking further.  We construct two 
right-hand sides, each coming from the span of some subset of eigenvectors.  We can 
further specify how many eigenvector components they have in common and see how
this affects convergence.

Let $\curl{\vek q_{1},\vek q_{2},\ldots, \vek q_{n}}$ be the orthonormal eigenvectors of $\vek A$, 
in ascending order according to the magnitude of their associated eigenvalues.  We define
the following subspaces,
\begin{eqnarray}
	\CQ_{S_{1}} & = & \span\curl{\vek q_{1},\ldots,\vek q_{100}}\nonumber\\
	\CQ_{S_{2}} & = & \span\curl{\vek q_{101},\ldots,\vek q_{200}}\nonumber\\
	\CQ_{L} & = &  \span\curl{\vek q_{n-200+1},\ldots,\vek q_{n}}\nonumber
\end{eqnarray}
 
\begin{figure}[htb]
\begin{center}

\includegraphics[scale=0.4]{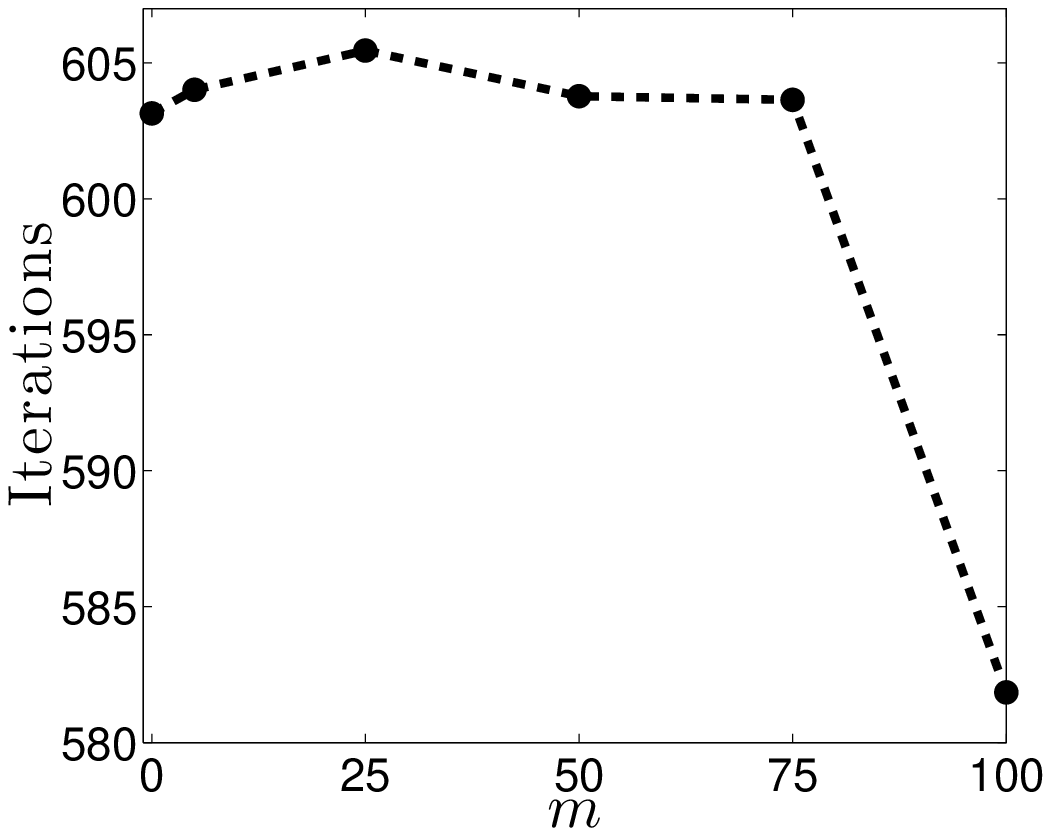}
\end{center}
\vspace*{-1em}
\caption{Average performance of our method for pairs of right-hand sides built from 
the $200$ eigenvectors associated with the smallest magnitude eigenvalues, as described 
in (\ref{eqn.smalleigsRHS}).  The right-hand sides have components from $2m$ common eigenvectors,
but are otherwise orthogonal.  For each $m$, the iteration counts and times are the averages of 
the results of tests for $100$ such pairs of right-hand sides.}
\label{fig.blockMINRES-eigmix-small}
\end{figure}

 In the first experiment, we construct both right-hand sides from eigenvectors associated
 only to eigenvalues of smaller magnitude, i.e., $\vek b_{1},\vek b_{2}\in\CQ_{S_{1}} \oplus\CQ_{S_{2}} $,
 such that a fixed number of eigenvectors are used to construct both vectors.  We define the two right-hand
 sides
 \begin{equation}\label{eqn.smalleigsRHS}
\vek b_{1} = \sum_{i=1}^{100+m}\alpha_{i}\vek q_{i} \mbox{\ \ and\ \ }\vek b_{2} = \sum_{i=100-m+1}^{200}\beta_{i}\vek q_{i}.
\end{equation}
For $m=0$, $\vek b_{1}$ and $\vek b_{2}$ are orthogonal.  For $m=1$, they both have components
from $\vek q_{100}$ and $\vek q_{101}$ but are otherwise orthogonal.  For $m=100$, both right-hand
sides have components in all $200$ basis vectors of $\CQ_{S_{1}} \oplus\CQ_{S_{2}}$.  For various
values of $m$, we can test the performance of our algorithm.  The coefficients $\curl{\alpha_{i}}$ and 
$\curl{\beta_{i}}$ are generated using Matlab's \verb|rand()| command. In order to avoid judging performance
based on a specific random example (which may be an outlier), for each $m$ tested, we generated $100$
different pairs of right-hand sides.  In Figure \ref{fig.blockMINRES-eigmix-small}, we plot the average iteration counts
over the $100$ tests for each $m$.  Until $m=100$, we see little change in the iteration counts.

\begin{figure}[htb]
\begin{center}
\includegraphics[scale=0.4]{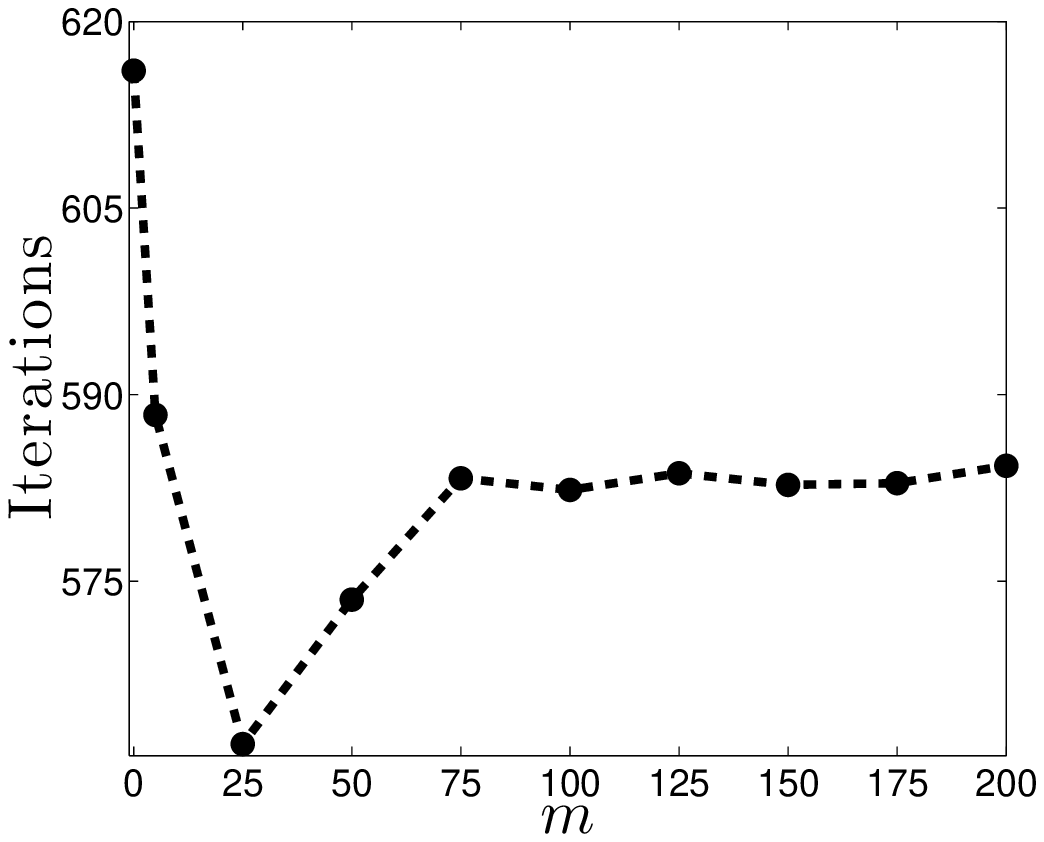}
\end{center}
\vspace*{-1em}
\caption{Experiments constructed in the same as those shown in Figure \ref{fig.blockMINRES-eigmix-small}
but now with pairs of right-hand sides constructed according to (\ref{eqn.largeeigsRHS}).  Again
the average iterations and times from tests for $100$ pairs of right-hand sides are given.}
\label{fig.blockMINRES-eigmix-bigsmall}
\end{figure}

We also performed the same experiment but constructed the two right-hand sides using
eigenvectors from different parts of the spectrum.  For different values of $m$, we
define
 \begin{equation}\label{eqn.largeeigsRHS}
\vek b_{1} = \sum_{i=1}^{200}\alpha_{i}\vek q_{i} + \underbrace{\sum_{i=n_{1}-200+1}^{n_{1}-200+m}\alpha_{i}\vek q_{i}}_{=0\mbox{\ when\ }m=0}\mbox{\ \ and\ \ }\vek b_{2} = \underbrace{\sum_{i=201-m}^{200}\beta_{i}\vek q_{i}}_{=0\mbox{\ when\ }m=0}+\sum_{i=n_{1}-200+1}^{n_{1}}\beta{i}\vek q_{i}.
\end{equation}
When $m=0$, we have $\vek b_{1}\in\CQ_{S_{1}} \oplus\CQ_{S_{2}}$ and 
$\vek b_{2}\in\CQ_{L}$ and they are orthogonal.  For $m=1$ they share components
from two eigenvectors ($\vek q_{200}$ and $\vek q_{n_{1}-200+1}$).  For $m=200$, both right-hand 
sides have components from every basis vector of $\CQ_{S_{1}} \oplus\CQ_{S_{2}}\oplus\CQ_{L}$.
As in the previous experiment, 100 random pairs of right-hand sides were generated for each $m$, and
the results averaged.  Average iterations counts are shown in 
Figure \ref{fig.blockMINRES-eigmix-bigsmall}.  We see a quick drop in iterations at $m=25$
followed by an increase.  Over all, mixing eigencomponents in this experiment produces a decrease
in iteration counts.  

This is by no means a rigorous analysis of the convergence of a block method.  These experiments only are meant 
to illustrate the variability of performance of a block method for different right-hand sides 
and provide some insight into this phenomenon.

\section{Conclusions}\label{section.conclusions}

We have presented an implementation of the block MINRES algorithm  
based on the band Lanczos process.  This version is designed to perform many
operations in a block fashion while maintaining the band Lanczos method's
easy-to-implement breakdown detection property.  
We provide not only a theoretical derivation of the algorithm but also a discussion of the 
practical implementation issues which need to be addressed to fully take advantage of the 
efficiencies which arise in a block method for symmetric systems.  This variant of the block 
MINRES method handles dependence of block Krylov subspace basis vectors in a more 
straightforward manner than a block Lanczos-based algorithm. A software implementation in Matlab is provided 
at \url{http://math.soodhalter.com/software.php}.

\section*{Acknowledgment}\label{section.acknowledgement}

The author would like to thank Sebastian Birk, Michael Parks, and Daniel Szyld  
for their constructive editorial comments and
suggestions.  The author would also like to express gratitude to the two reviewers and editor who
offered extensive comments and constructive criticism, which were of great help 
in improving this manuscript.  In particular, it should be noted that the second reviewer
suggested the expression "dynamic substitutes bench" to describe the random vectors
stored for use in maintaining the block size.

\bibliographystyle{siam}

\begin{thebibliography}{10}

\bibitem{ABFH.2000}
{\sc Jos\'{e}~I. Aliaga, Daniel~L. Boley, Roland~W. Freund, and Vicente
  Hern\'{a}ndez}, {\em A {L}anczos-type method for multiple starting vectors},
  Mathematics of computation, 69 (2000), pp.~1577--1602.

\bibitem{B.2000}
{\sc James Baglama}, {\em Dealing with linear dependence during the iterations
  of the restarted block {L}anczos methods}, Numerical Algorithms, 25 (2000),
  pp.~23--36.

\bibitem{BDJ.2006}
{\sc Allison~H. Baker, John~M. Dennis, and Elisabeth~R. Jessup}, {\em On
  improving linear solver performance: a block variant of {GMRES}}, SIAM J.
  Sci. Comput., 27 (2006), pp.~1608--1626.

\bibitem{BF.2013}
{\sc Sebastian Birk and Andreas Frommer}, {\em A deflated conjugategate
  gradient method for multiple right-hand sides and multiple shifts},  (In
  preparation).

\bibitem{DS.2012}
{\sc J.~Jack Dongarra and Aad~J. van~der Sten}, {\em High-performance computing
  systems: status and outlook}, Acta Numer., 21 (2012), pp.~379--474.

\bibitem{D.2001}
{\sc Augustin~A. Dubrulle}, {\em Retooling the method of block conjugate
  gradients}, Electronic Transactions on Numerical Analysis, 12 (2001),
  pp.~216--233 (electronic).

\bibitem{F.1995}
{\sc Roland~W. Freund}, {\em Computation of matrix {P}ad\'e approximations of
  transfer functions via a {L}anczos-type process}, in Approximation theory
  {VIII}, {V}ol.\ 1 ({C}ollege {S}tation, {TX}, 1995), vol.~6 of Ser. Approx.
  Decompos., World Sci. Publ., River Edge, NJ, 1995, pp.~215--222.

\bibitem{FM.1997}
{\sc Roland~W. Freund and Manish Malhotra}, {\em A block {QMR} algorithm for
  non-{H}ermitian linear systems with multiple right-hand sides}, in
  Proceedings of the {F}ifth {C}onference of the {I}nternational {L}inear
  {A}lgebra {S}ociety ({A}tlanta, {GA}, 1995), vol.~254, 1997, pp.~119--157.

\bibitem{freund.QMR.1991}
{\sc Roland~W. Freund and No{\"e}l~M. Nachtigal}, {\em Q{MR}: a quasi-minimal
  residual method for non-{H}ermitian linear systems}, Numerische Mathematik,
  60 (1991), pp.~315--339.

\bibitem{Greenbaum1997a}
{\sc Anne Greenbaum}, {\em Iterative Methods for Solving Linear Systems}, SIAM,
  Philadelphia, 1997.

\bibitem{Gutknecht2008}
{\sc Martin~H. Gutknecht and Thomas Schmelzer}, {\em Updating the {QR}
  decomposition of block tridiagonal and block {H}essenberg matrices}, Applied
  Numerical Mathematics, 58 (2008), pp.~871--883.

\bibitem{GS.2009}
\leavevmode\vrule height 2pt depth -1.6pt width 23pt, {\em The block grade of a
  block {K}rylov space}, Linear Algebra and its Applications, 430 (2009),
  pp.~174--185.

\bibitem{H.2010}
{\sc Mark Hoemmen}, {\em Communication-avoiding Krylov subspace methods}, PhD
  thesis, University of California Berkeley, 2010.

\bibitem{L.2003}
{\sc Julian Langou}, {\em Iterative methods for solving linear systems with
  multiple right-hand sides}, PhD thesis, CERFACS, France, 2003.

\bibitem{L.2006}
{\sc Damian Loher}, {\em Reliable nonsymmetric block {L}anczos algorithms}, PhD
  thesis, Diss.\ no.\ 16337, ETH Zurich, Zurich, Switzerland, 2006.

\bibitem{OLeary1980}
{\sc Dianne~P. O'Leary}, {\em {The block conjugate gradient algorithm and
  related methods}}, Linear Algebra and its Applications, 29 (1980),
  pp.~293--322.

\bibitem{Paige1975}
{\sc Chris~C. Paige and Michael~A. Saunders}, {\em Solutions of sparse
  indefinite systems of linear equations}, SIAM Journal on Numerical Analysis,
  12 (1975), pp.~617--629.

\bibitem{PSS.2013}
{\sc Michael~L. Parks, Kirk~M. Soodhalter, and Daniel~B. Szyld}, {\em Block
  krylov subspace recycling}, In Preparation.

\bibitem{RS.2006}
{\sc Micka{\"e}l Robb{\'e} and Miloud Sadkane}, {\em Exact and inexact
  breakdowns in the block {GMRES} method}, Linear Algebra and its Applications,
  419 (2006), pp.~265--285.

\bibitem{Ruhe1979}
{\sc Axel Ruhe}, {\em Implementation aspects of band {L}anczos algorithms for
  computation of eigenvalues of large sparse symmetric matrices}, Mathematics
  of Computation, 33 (1979), pp.~680--687.

\bibitem{Saad.Iter.Meth.Sparse.2003}
{\sc Yousef Saad}, {\em Iterative Methods for Sparse Linear Systems}, SIAM,
  Philadelphia, {S}econd~ed., 2003.

\bibitem{Saad.GMRES.1986}
{\sc Yousef Saad and Martin~H. Schultz}, {\em {GMRES}: A generalized minimal
  residual algorithm for solving nonsymmetric linear systems}, SIAM Journal on
  Scientific and Statistical Computing, 7 (1986), pp.~856--869.

\bibitem{S.2004}
{\sc Thomas Schmelzer}, {\em Block {Krylov} methods for {Hermitian} linear
  systems}, master's thesis, 2004.

\bibitem{Simoncini.Conv-Block-GMRES}
{\sc Valeria Simoncini and Efstratios Gallopoulos}, {\em Convergence properties
  of block {GMRES} and matrix polynomials}, Linear Algebra and its
  Applications, 247 (1996), pp.~97--119.

\bibitem{szyld.simoncini.survey.2007}
{\sc Valeria Simoncini and Daniel~B. Szyld}, {\em Recent computational
  developments in {K}rylov subspace methods for linear systems}, Numerical
  Linear Algebra with Applications, 14 (2007), pp.~1--59.

\bibitem{S.2014-arxiv}
{\sc Kirk~M. Soodhalter}, {\em A block {MINRES} algorithm based on the band
  lanczos method}, Tech. Report 1301.2102v2, arXiv, 2013.

\bibitem{Vital1990}
{\sc Brigitte Vital}, {\em Etude de quelques m{\'e}thodes de r{\'e}solution de
  probl\`{e}mes lin{\'e}aires de grande taille sur multiprocesseur}, PhD
  thesis, Universit{\'e} de Rennes, 1990.

\end{thebibliography}

\end{document}